\documentclass[12pt]{amsart}
\usepackage{epsf}
\usepackage{amsmath}
\usepackage{amssymb}
\usepackage{epic}
\usepackage{graphics}
\newtheorem{theo}{\bf Theorem}

\newtheorem{conj}[theo]{\bf Conjecture}
\newtheorem{defi}[theo]{\bf Definition}

\newcommand{\twelvedots}{
  \put(0,0){\circle*{4}}
  \put(0,20){\circle*{4}}
  \put(0,40){\circle*{4}}

  \put(40,0){\circle*{4}}
  \put(40,20){\circle*{4}}
  \put(40,40){\circle*{4}}

  \put(80,0){\circle*{4}}
  \put(80,20){\circle*{4}}
  \put(80,40){\circle*{4}}

  \put(120,0){\circle*{4}}
  \put(120,20){\circle*{4}}
  \put(120,40){\circle*{4}}
}
\newcommand{\sixteendots}{
  \put(0,0){\circle*{4}}
  \put(0,20){\circle*{4}}
  \put(0,40){\circle*{4}}
  \put(0,60){\circle*{4}}

  \put(40,0){\circle*{4}}
  \put(40,20){\circle*{4}}
  \put(40,40){\circle*{4}}
  \put(40,60){\circle*{4}}

  \put(80,0){\circle*{4}}
  \put(80,20){\circle*{4}}
  \put(80,40){\circle*{4}}
  \put(80,60){\circle*{4}}

  \put(120,0){\circle*{4}}
  \put(120,20){\circle*{4}}
  \put(120,40){\circle*{4}}
  \put(120,60){\circle*{4}}
}

\begin{document}

\title{Conjectures on three-dimensional stable matching}
\author{Kimmo Eriksson}
\address{IMA, M{\"a}lardalens h{\"o}gskola \\
   Box 883 \\
   SE-721 23 V{\"a}ster{\aa}s, Sweden}
\email{kimmo.eriksson@mdh.se}
\author{Jonas Sj{\"o}strand}
\address{Dept.~of Mathematics \\
   KTH \\
   SE-100 44 Stockholm, Sweden}
\email{jonass@kth.se}
\author{Pontus Strimling}
\address{IMA, M{\"a}lardalens h{\"o}gskola \\
   Box 883 \\
   SE-721 23 V{\"a}ster{\aa}s, Sweden}
\email{pontus.strimling@mdh.se} 
\date{October 31, 2004}

\begin{abstract}
We consider stable three-dimensional matchings of three categories of agents, such as women, men and dogs.
This was suggested long ago by Knuth (1976), but very little seems to have been published on this problem.
Based on computer experiments, we present a couple of conjectures as well as a few counter-examples to other
natural but discarded conjectures. In particular, a \emph{circular} 3D matching is one where women only care
about the man, men only care about the dog, and dogs only care about the woman they are matched with. We
conjecture that a stable outcome always exists for any circular 3D matching market, and we prove it for
markets with at most four agents of each category.
\end{abstract}
\maketitle

\section{Introduction}\noindent
The \emph{stable marriage problem} is a well-known problem in matching theory: Given a set of men and a set of women,
find a matching that is stable. A matching is stable if there is no blocking pair, that is, a man and a woman who
would prefer each other to their current partners in the matching. Gale and Shapley (1962) introduced this problem
and proved that a stable matching must exist by describing an algorithm that produces such a matching.

The theory of stable matchings has become an important subfield within game theory, as documented by the
book of Roth and Sotomayor (1990). Game theorists are interested in the applications of matching theory in real
markets. However, the theory of stable matchings also appeal to combinatorialists and computer scientists. Indeed,
the first book on the subject was written by combinatorial computer scientist \emph{extraordinaire} Donald E.~Knuth
(Knuth, 1976). These books discuss not only the two-sided matching problem, but also the one-sided so called
\emph{roommate problem} where any two agents can form a pair. For the roommate problem, it is easy to find
counter-examples to stability.

Knuth (1976) lists a dozen suggested further directions for
research on stable matchings, one of which is to investigate
three-dimensional stable matching, say of women, men and dogs. A
3D matching is a partition of the agents into triples consisting
of one agent of each type. A 3D matching is stable if no blocking
triple exists, that is, a triple $(a,b,c)$ which each of its
members would prefer to their current triples in the matching. Of
course, this calls for a definition of how agents rank triples
based on their preferences on other individuals. Many
possibilities exist.

The only paper we have found on three-dimensional stable matching
problems is a complexity investigation by Ng and Hirschberg
(1991). They showed that some instances of what they dub the
\emph{three-gender stable marriage problem} (3GSM) do not have any
strongly stable outcomes, and proved that the decision problem is
NP-complete. (Strong stability is defined in the next section.) As
an open problem, Ng and Hirschberg mention the \emph{circular}
3GSM where, say, women rank triples based only on the man in the
triple, and similarly men care only about dogs, and dogs care only
about women. The origin of this problem is attributed to Knuth.

In this paper, we will report the results of our investigation of the existence of stable outcomes of 3GSM in the
circular version as well as in a few other versions. Our main conjecture is that circular 3GSM always allows a
stable outcome. We prove this for all instances with at most four agents of each gender, and we describe the
evidence from computer experimentation that led us to this conjecture.

\section{Problem definition}\noindent
Let $N$ be the maximal number of agents of each of the three genders. Thus, $N=3$ means that we have at most
nine agents (three women, three men, three dogs). Without loss of generality we can assume that we have
the maximal number of agents of each gender, for otherwise we can just fill the ranks with dummy agents who
everybody likes less than any real agent.

We will assume that every agent ranks all possible triples based
on her ranking of the individuals of the other genders. For
example, every woman will have a preference list of length $2N$ on
the set of all men and dogs. If $a$ prefers $b_1$ to $c_1$ to
$b_2$, we write this as $b_1 >_a c_1 >_a b_2$. Agents together
with their preferences constitute a \emph{3G matching market}.

We will see several different rules for generating preferences on
triples from preferences on individuals. If $a$ prefers triple
$T_1$ to $T_2$ we write $T_1 >_a T_2$. If preferences are not
strict, we use $\ge_a$ instead (\emph{weak} preference).

A \emph{3G matching} is a partition of the agents into $N$ triples consisting of a woman, a man and a dog.
Given agents' preferences on triples, a 3G matching has a \emph{blocking triple} $T$ if all members of $T$
strictly prefer $T$ to their current triples in the matching. A 3G matching is \emph{stable} if it has no
blocking triples. A triple $T$ is \emph{weakly blocking} if some member strictly prefers $T$ to her current
matching and the other members weakly prefer $T$. A 3G matching is \emph{strongly stable} if it has no weakly
blocking pair.

The problem is: Given a triple preference rule, does a stable 3G
matching exist for every matching market of size $N$? We have been
interested in the following rules for preferences on triples. Let
$T=(a,b,c)$ and $T'=(a,b',c')$.
\begin{itemize}
\item \emph{Circular}: $T >_a T'$ if and only if $b >_a b'$. Similarly, men care only about dogs, and dogs care
only about women.
\item \emph{Weakest link}: $T >_a T'$ if $\min_a(b,c) >_a \min_a(b',c')$. In general, agents rank triples according
to the weakest link of the triple, that is, according to their
least preferred partner.
\item \emph{Strongest link}: $T >_a T'$ if $\max_a(b,c) >_a \max_a(b',c')$. In general, agents rank triples according
according to their most preferred partner.
\end{itemize}
Note that the circular rule is a special case both of the weakest link rule and the strongest link rule, depending
on whether we let the gender cared about in the circular rule be consistently low-ranked or high-ranked
in relation to the other gender.

\section{Investigating stability by computer}\noindent
For size $N$, the number of different matching markets is $(2N)!^{3N}$, since each agent ranks
all $2N$ agents of other genders. Even if isomorphic copies were deleted, the number of markets is
large already for $N=3$ and daunting for $N > 3$.

In order to investigate stability, we wrote a program (in Python) which starts by generating a random
market of a given size $N\le 5$. For this market each of the $N!^2$ possible matchings are checked for
stability (according to a given triple preference rule).
The number of stable matchings is recorded. Then a local search for markets with fewer stable
matchings is carried out as follows: Each of the $3N$ agents, in turn, changes its preference list to
every possible alternative permutation. For each of these markets the number of stable matchings is computed,
and whenever a new minimum is found, search proceeds from this market.

If a market with zero stable matchings is found, we have a counter-example to the existence of stable matchings
for the given preferences rule.

If no market with zero stable matchings is found, we have an indication that there is none. More
specifically we obtain an indication of the minimum number of stable matchings.

We started searching with $N=3$, in which case we never found any markets without stable matchings for any
of the three preference rules. This came as a surprise to us, and so we continued with $N=4$.

\section{Counter-examples to stability}\noindent
With $N=4$, our computer search found a \emph{counter-example to
stability under the weakest link rule}. In order to present the
preferences in a convenient way, we write lists of ranking numbers
from 1 to 8. Lower rank means more preferred.
\begin{align*}
a_1 &: [[2, 8, 3, 7], [5, 1, 4, 6]] \\
a_2 &: [[7, 4, 2, 1], [3, 8, 5, 6]] \\
a_3 &: [[8, 4, 6, 3], [5, 2, 1, 7]] \\
a_4 &: [[1, 2, 8, 6], [7, 5, 3, 4]] \\
&\\
b_1 &: [[6, 5, 7, 2], [8, 3, 1, 4]] \\
b_2 &: [[7, 1, 4, 6], [5, 3, 8, 2]] \\
b_3 &: [[4, 3, 1, 2], [7, 6, 5, 8]] \\
b_4 &: [[8, 7, 5, 4], [6, 2, 1, 3]] \\
&\\
c_1 &: [[8, 3, 2, 7], [5, 4, 6, 1]] \\
c_2 &: [[5, 1, 6, 4], [8, 7, 2, 3]] \\
c_3 &: [[8, 1, 3, 2], [5, 6, 4, 7]] \\
c_4 &: [[8, 7, 5, 3], [1, 6, 2, 4]]
\end{align*}
For women, the lists give the ranking of $[[b_1, b_2, b_3, b_4],
[c_1, c_2, c_3, c_4]]$. For men, the lists denote the ranking of
$[[a_1, a_2, a_3, a_4], [c_1, c_2, c_3, c_4]]$. For dogs, it is
the ranking of $[[a_1, a_2, a_3, a_4], [b_1, b_2, b_3, b_4]]$.

Furthermore, we found a \emph{counter-example to strong stability under the circular rule}:
\begin{align*}
a_1 &: [1, 2, 4, 3] \\
a_2 &: [2, 1, 3, 4] \\
a_3 &: [2, 4, 3, 1] \\
a_4 &: [1, 3, 4, 2] \\
&\\
b_1 &: [3, 1, 2, 4] \\
b_2 &: [4, 1, 2, 3] \\
b_3 &: [3, 2, 4, 1] \\
b_4 &: [2, 4, 3, 1] \\
&\\
c_1 &: [2, 3, 4, 1] \\
c_2 &: [1, 4, 2, 3] \\
c_3 &: [2, 3, 4, 1] \\
c_4 &: [3, 1, 2, 4]
\end{align*}

\section{Conjectures}\noindent
In the light of the counter-examples of the previous section, it is intriguing that days of computer time
for $N=4$ and $N=5$ have not resulted in any counter-examples to the following conjecture.
\begin{conj}
Under the strongest link rule, every 3G matching market has a stable matching.
\end{conj}\noindent
In fact, since the lowest number of stable matchings found by the computer is \emph{two} for both cases
$N=4$ and $N=5$ (examples of such markets available from the authors),
we propose the following stronger conjecture.
\begin{conj}
Under the strongest link rule, every 3G matching market has at least two stable matchings.
\end{conj}\noindent
For the circular rule, which is a special case of the strongest link rule, the computer always find many
stable matchings so we guess that the minimal number of stable matchings increase with $N$. We do not
have enough evidence to produce a firmer guess, so our main conjecture is simply:
\begin{conj}
Under the circular rule, every 3G matching market has a stable matching.
\end{conj}\noindent
The circular case conjecture would seem to be amenable to an algorithmic approach similar to the deferred
acceptance procedure of Gale and Shapley (1962) for two-sided marriages. In other words, it is not difficult
to come up with ideas of algorithms like the following:

``Let all women propose to the men they prefer most. Let every man
tentatively accept the woman who is most preferred by the dog the man prefers most. Continue until all women
propose to different men.''

However, all such ideas seem to run into problems and we have
resigned ourselves to nonconstructive approaches: In the next
section we will describe how one might try to apply Scarf's
theorem (Scarf, 1971) on balanced games, and why this fails.
Finally we will carry out a case-by-case analysis, which is doable
for $N\le 4$ but then seems to get out of hand.

\section{Circular 3GSM is not a balanced game}\noindent
If you want to show that a game has a nonempty core (which in our
case is equivalent to the existence of a stable matching), one
approach is to show that it is {\em balanced} in the sense of
Scarf~\cite{scarf}.
\begin{theo}[{\bf Scarf, 1967}]
A balanced $n$-person game always has a nonempty core.
\end{theo}\noindent
Quinzii~\cite{quinzii} showed (in a more general setting) that the
usual two-dimensional matching game is balanced. However, in this
section we will see that our three-dimensional matching game is
not.

\begin{figure}
\begin{center}
\setlength{\unitlength}{0.5mm}
\begin{picture}(100,100)(0,0)
\put(0,0){\resizebox{50mm}{!}{\includegraphics{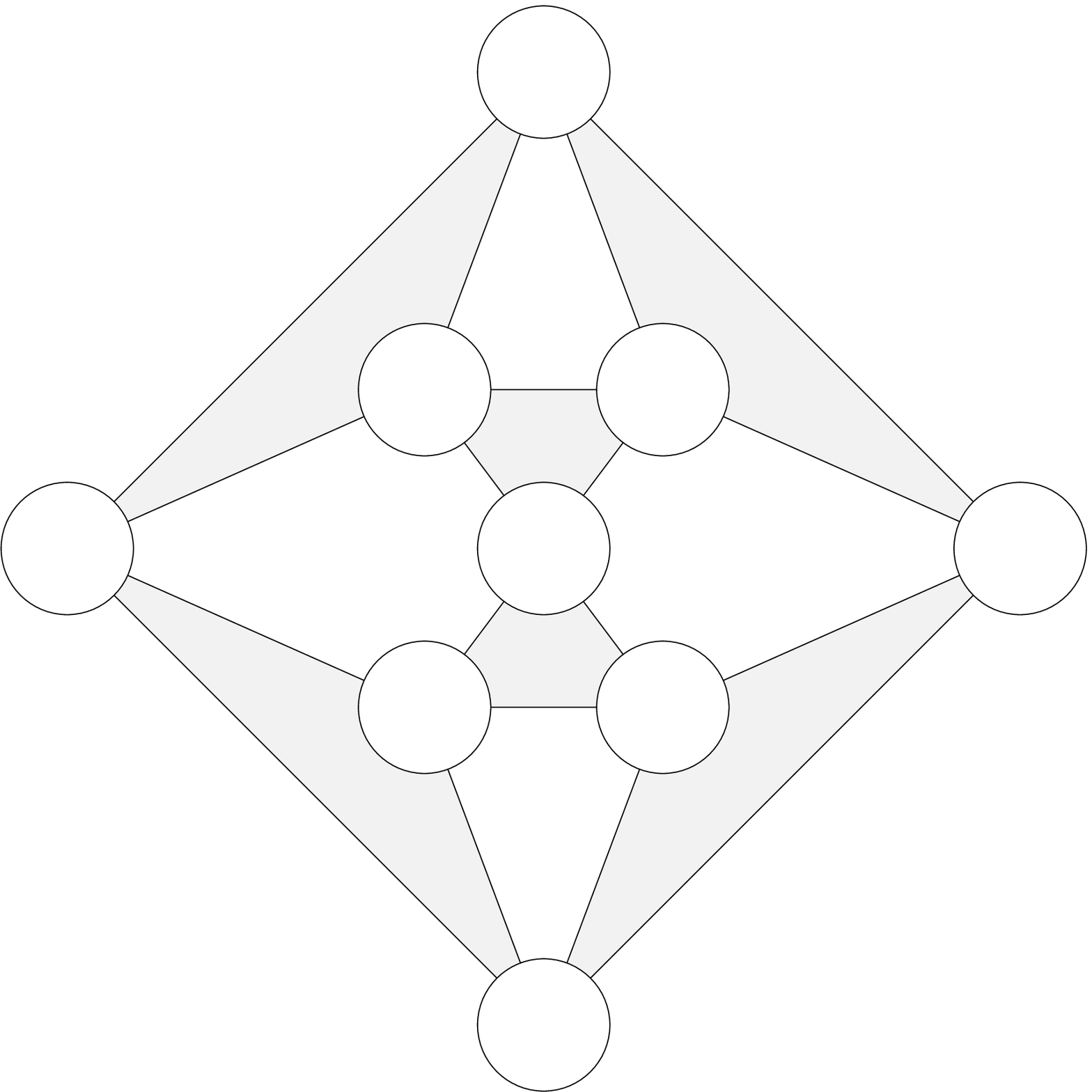}}}
\put(46.5,4){$a_1$} \put(36,33){$b_1$} \put(58,33.3){$c_1$}
\put(46.5,48){$a_2$} \put(36,62){$b_2$} \put(58,62.5){$c_2$}
\put(46.5,92){$a_3$} \put(3,48){$c_3$} \put(91,47.5){$b_3$}
\end{picture}
\caption{An example showing that our game is not balanced.}
\label{fig:counterexample}
\end{center}
\end{figure}

For the general definition of a balanced game we refer
to~\cite{scarf}. Here, we will merely examine what it would mean
for our game to be balanced.
\begin{defi}
A collection $C$ of triples is {\em balanced} if there is possible
to find nonnegative real weights $\delta_T$, for each triple $T$
in $C$, such that, for each person $x$,
$$\sum_{x\in T\in C}\delta_T=1.$$
\end{defi}
\noindent A {\em utility vector} is a list where every person has
written down her utility goal, that is, how happy she hopes to
become. A utility vector is {\em realizable} if there is a
matching such that every person reach her utility goal. A utility
vector is {\em realizable for a triple} if all people in the
triple would reach their utility goal if the triple were formed.
\begin{defi}
Our game is {\em balanced} if, for every balanced collection $C$
of triples, a utility vector is realizable if it is realizable for
every triple in $C$.
\end{defi}
\noindent Now, we present a counterexample of size $3+3+3$. Let
$C$ be the collection of triples corresponding to the shaded
triangles in figure~\ref{fig:counterexample}.
This collection is balanced, since every person belongs to exactly
two triples (let all $\delta_T=1/2$). Choose the preferences so
that the edges in the figure correspond to rank 1 or 2. For
example, $a_1$ will rank $b_1$ and $b_3$ as number 1 and 2 (in any
order) and $b_2$ as number 3. Now consider the utility vector
where every person hopes to get at least her second-best choice.
This is obviously realizable for every triple in $C$, so if the
game were balanced, the utility vector would be realizable.
Since every instance
of ``$x$ ranks $y$ as number 1 or 2'' has a corresponding edge in
the figure, a realization of the utility vector is equivalent to a
disjoint family of triangles (not necessarily shaded) in the
figure which covers all people. But there is no such family: To
cover $a_2$, either of the triangles $a_2b_1c_1$ and $a_2b_2c_2$
must belong to the family. But none of the three triangles
containing $a_1$ is disjoint with $a_2b_1c_1$, and none of the
three triangles containing $a_3$ is disjoint with $a_2b_2c_2$.

\section{Proof of stability for circular 3GSM with $N\le 4$}

\subsection{Notation} \noindent Say that a person is {\em
i-content} if she has got the person she ranks as number $i$.

If $x$'s favorite $y$'s favorite $z$ ranks $x$ as number $i$, we
say that $xyz$ is a $11i$-triple.
The proposition ``there is no $11i$-triple for $i<j$'' is called
the {\em $11j$-condition}.

We will often use a dot diagram to describe partial information of
the preferences. Here is an example:
\begin{center}
  \begin{picture}(130,50)(-4,-7)
    \twelvedots
    \put(0,40){\line(1,0){80}}
    \put(80,40){\line(2,-1){40}}
    \dashline[-7]{3}(0,40)(40,20)
    \put(99,31){\tiny 3}

    \put(-4,-7){\tiny $a_3$}
    \put(-4,13){\tiny $a_2$}
    \put(-4,33){\tiny $a_1$}
    \put(36,-8){\tiny $b_3$}
    \put(36,12){\tiny $b_2$}
    \put(36,32){\tiny $b_1$}
    \put(76,-7){\tiny $c_3$}
    \put(76,13){\tiny $c_2$}
    \put(76,33){\tiny $c_1$}
    \put(116,-7){\tiny $a_3$}
    \put(116,13){\tiny $a_2$}
    \put(116,33){\tiny $a_1$}
  \end{picture}
\end{center}
A solid line means ``rank 1'', a dashed line means ``rank 2'' and
an $i$-labeled line means ``rank $i$''. So from the diagram above
we get the information that $a_1$'s second-best choice is $b_2$,
but her favorite is $b_1$ whose favorite is $c_1$ who ranks $a_2$
as number 3. To make the following pages more readable, we will
omit the dot labels, always implicitly referring to the labeling
above (in the 4-case, of course, there will be an additional row
$a_4b_4c_4a_4$ at the bottom).

\subsection{The 3-case}\noindent
\begin{theo}\label{thm:3}
In the 3-case there is always a stable matching. Furthermore, for
any person $x$, there is a stable matching such that either $x$ is
1-content, or $x$'s favorite is 1-content and $x$ is 2-content.
\end{theo}
\begin{proof}
We can assume that $x=a_1$ whose favorite is $b_1$ whose favorite
is $c_1$.

Suppose there is a 111-triple $abc$. Then pick that triple to the
matching. Since $a$, $b$ and $c$ all are 1-content we can choose
the other two triples however we like --- the resulting matching
will be stable anyway. If $b_1=b$ he is 1-content and we choose
the remaining two triples so that $a_1$ gets her first or second
choice. If $b_1\neq b$ then $a_1\neq a$ and we let the first
triple be $a_1b_1C$ where $C$ is $b_1$'s first or second choice,
and the second triple be the three remaining people.

In the following we assume there is no 111-triple, that is,
we assume the 112-condition. Say, without loss of generality, that
$c_1$'s favorite is $a_2$.
\begin{center}
  \begin{picture}(130,50)(-4,-4)
    \twelvedots
    \put(0,40){\line(1,0){80}}
    \put(80,40){\line(2,-1){40}}
  \end{picture}
\end{center}
By the 112-condition, $a_2$'s favorite is not $b_1$, so we can
assume that it is $b_2$. Now, $b_2$'s favorite is not $c_1$, so we
can assume that it is $c_2$.
\begin{center}
  \begin{picture}(130,50)(-4,-4)
    \twelvedots
    \put(0,40){\line(1,0){80}}
    \put(80,40){\line(2,-1){40}}
    \put(0,20){\line(1,0){80}}
  \end{picture}
\end{center}
Now the matching $a_1b_1c_1,a_2b_2c_2,a_3b_3c_3$ is stable: $a_1$,
$a_2$, $b_1$ and $b_2$ are 1-content, so a blocking triple must
contain both $a_3$ and $b_3$ which already are together.
\end{proof}

\subsection{The 4-case}\noindent
\begin{theo}
In the 4-case there is always a stable matching.
\end{theo}
\noindent The proof is a technical case study who will last for
the rest of this section.

We have the following cases:
\begin{description}
\item[The 111-case] There is a 111-triple.
\item[The 112-case] There is no 111-triple, but there is a 112-triple.
\item[The 113-case] There is no 111- or 112-triple, but there is
  a 113-triple.
\item[The 114-case] There is no 111-, 112- or 113-triple.
\end{description}
\noindent The 111-case is trivial: remove the 111-triple and find
a stable 3-matching of the remainder; then the 111-triple together
with the 3-matching is a stable 4-matching since the people in the
111-triple are 1-content.

The 114-case is also simple: If there is a person $x$ who is the
favorite of at least two people, then $x$'s favorite must rank
some of these people as number 1, 2 or 3, and we get a 111-, 112-,
or 113-triple. Thus, in the 114-case we know that no two people
have the same favorite. We just let all $a_i$ be 1-content and we
get a stable matching.

\subsubsection{The 112-case} We have the following situation.
\begin{center}
  \begin{picture}(130,70)(-4,-4)
    \sixteendots
    \put(0,60){\line(1,0){80}}
    \put(80,60){\line(2,-1){40}}
    \dashline[-7]{3}(80,60)(120,60)
  \end{picture}
\end{center}
By the 112-condition, $a_2$'s favorite is not $b_1$, so we can
assume it is $b_2$. Now $b_2$'s favorite is not $c_1$, so we can
assume it is $c_2$.
\begin{center}
  \begin{picture}(130,70)(-4,-4)
    \sixteendots
    \put(0,60){\line(1,0){80}}
    \put(80,60){\line(2,-1){40}}
    \dashline[-7]{3}(80,60)(120,60)
    \put(0,40){\line(1,0){80}}
  \end{picture}
\end{center}
We remove the triple $a_1b_1c_1$ for a while. By
theorem~\ref{thm:3} there is a stable 3-matching of the remaining
people such that either $a_2$ is 1-content, or $b_2$ is 1-content
and $a_2$ is 2-content. This 3-matching forms a 4-matching
together with the triple $a_1b_1c_1$. We will show that this
4-matching is stable.

Suppose there is a blocking triple. It has to contain someone
among $a_1$, $b_1$ and $c_1$. Since $a_1$ and $b_1$ are 1-content
they do not belong to the blocking triple, so $c_1$ does. The only
person $c_1$ wants to switch to is $a_2$, so $a_2$ belongs to the
blocking triple. Then $a_2$ cannot have her favorite $b_2$, so, by
construction of the 3-matching, $a_2$ has her favorite among $b_3$
and $b_4$, and $b_2$ has $c_2$. So $a_2$ wants to switch only to
$b_2$ or possibly $b_1$, both of which are 1-content.

\subsubsection{The 113-case} If some person $x$ is the favorite of
at least three people, then $x$'s favorite must rank some of these
people as number 1 or 2, and we get a 111- or a 112-triple.
Therefore, we can split the 113-case into two subcases:
\begin{description}
\item[Subcase 1] Every person is the favorite of either
  zero or two people.
\item[Subcase 2] There is a person who is the favorite of exactly
  one person, but every person is the favorite of at most two people.
\end{description}

\paragraph{\bf Subcase 1} Here we suppose that every person is
the favorite of either zero or two people. Then we have the
following situation.
\begin{center}
  \begin{picture}(130,70)(-4,-4)
    \sixteendots
    \put(0,0){\line(1,0){40}}
    \put(0,20){\line(2,-1){40}}
    \put(0,40){\line(2,1){40}}
    \put(0,60){\line(1,0){80}}
  \end{picture}
\end{center}
By the 113-condition we know that $c_1$ ranks $a_1$ and $a_2$ as
number 3 and 4, so $a_3$ and $a_4$ must be number 1 and 2.
\begin{center}
  \begin{picture}(130,70)(-4,-4)
    \sixteendots
    \put(0,0){\line(1,0){40}}
    \put(0,20){\line(2,-1){40}}
    \put(0,40){\line(2,1){40}}
    \put(0,60){\line(1,0){80}}
    \put(80,60){\line(1,-1){40}}
    \dashline[-7]{3}(80,60)(120,0)
  \end{picture}
\end{center}
We see that $b_4$'s favorite is not $c_1$, so we can assume it is
$c_4$. Using the 113-condition again, we obtain the following.
\begin{center}
  \begin{picture}(130,70)(-4,-4)
    \sixteendots
    \put(0,0){\line(1,0){80}}
    \put(0,20){\line(2,-1){40}}
    \put(0,40){\line(2,1){40}}
    \put(0,60){\line(1,0){80}}
    \put(80,60){\line(1,-1){40}}
    \put(80,0){\line(1,1){40}}
    \dashline[-7]{3}(80,60)(120,0)
    \dashline[-7]{3}(80,0)(120,60)
  \end{picture}
\end{center}
Now we use that every favorite person is the favorite of exactly
two people.
\begin{center}
  \begin{picture}(130,70)(-4,-4)
    \sixteendots
    \put(0,0){\line(1,0){80}}
    \put(0,20){\line(2,-1){40}}
    \put(0,40){\line(2,1){40}}
    \put(0,60){\line(1,0){80}}
    \put(80,60){\line(1,-1){40}}
    \put(80,0){\line(1,1){40}}
    \put(40,40){\line(2,1){40}}
    \put(40,20){\line(2,-1){40}}
    \put(80,20){\line(1,0){40}}
    \put(80,40){\line(1,0){40}}
    \dashline[-7]{3}(80,60)(120,0)
    \dashline[-7]{3}(80,0)(120,60)
  \end{picture}
\end{center}
$a_2$'s second-best choice cannot be $b_2$ or $b_3$ since that
would violate the 113-condition, so it must be $b_2$. By applying
the same reasoning to $a_3$, $b_1$ and $b_4$, we obtain the following
diagram.
\begin{center}
  \begin{picture}(130,70)(-4,-4)
    \sixteendots
    \put(0,0){\line(1,0){80}}
    \put(0,20){\line(2,-1){40}}
    \put(0,40){\line(2,1){40}}
    \put(0,60){\line(1,0){80}}
    \put(80,60){\line(1,-1){40}}
    \put(80,0){\line(1,1){40}}
    \put(40,40){\line(2,1){40}}
    \put(40,20){\line(2,-1){40}}
    \put(80,20){\line(1,0){40}}
    \put(80,40){\line(1,0){40}}
    \dashline[-7]{3}(80,60)(120,0)
    \dashline[-7]{3}(80,0)(120,60)
    \dashline[-7]{3}(0,20)(40,20)
    \dashline[-7]{3}(0,40)(40,40)
    \dashline[-7]{3}(40,0)(80,40)
    \dashline[-7]{3}(40,60)(80,20)
  \end{picture}
\end{center}
Now the matching $a_1b_1c_3,a_2b_2c_1,a_3b_3c_4,a_4b_4c_2$
is stable: $a_1$ and $a_4$ are
1-content, so a blocking triple must contain $a_2$ or $a_3$, say
$a_2$. But $a_2$ wants to switch only to $b_1$ who wants to switch
only to $c_1$ who already has got $a_2$. The same reasoning works
for $a_3$.

\vspace{3mm}
\paragraph{\bf Subcase 2} Here we suppose that
there is a person, say $b_3$, who is the favorite of exactly one
person, but every person is the favorite of at most two people.

If every $b_i$ is the favorite of exactly one person, then it is
trivial to find a stable matching: Just let all $a_i$ be
1-content. So we assume there is a $b_i$, say $b_1$, who is the
favorite of exactly two people. Then we have the following
situation:
\begin{center}
  \begin{picture}(130,70)(-4,-4)
    \sixteendots
    \put(0,0){\line(1,0){40}}
    \put(0,20){\line(1,0){40}}
    \put(0,40){\line(2,1){40}}
    \put(0,60){\line(1,0){40}}
  \end{picture}
\end{center}
By the 113-condition, $b_1$'s favorite must rank $a_1$ and $a_2$
as 3 and 4.
\begin{center}
  \begin{picture}(130,80)(-4,-4)
    \sixteendots
    \put(0,0){\line(1,0){40}}
    \put(0,20){\line(1,0){40}}
    \put(0,40){\line(2,1){40}}
    \put(0,60){\line(1,0){120}}
    \put(80,60){\line(2,-1){40}}
    \put(99,61){\tiny 3}
    \put(99,51){\tiny 4}
  \end{picture}
\end{center}
Now, $b_3$ or $b_4$ cannot have $c_1$ as a favorite, since $c_1$
ranks $a_3$ and $a_4$ as 1 and 2, and that would violate the
113-condition. Thus we can assume that $b_3$'s favorite is, say,
$c_3$.
\begin{center}
  \begin{picture}(130,80)(-4,-4)
    \sixteendots
    \put(0,0){\line(1,0){40}}
    \put(0,20){\line(1,0){80}}
    \put(0,40){\line(2,1){40}}
    \put(0,60){\line(1,0){120}}
    \put(80,60){\line(2,-1){40}}
    \put(99,61){\tiny 3}
    \put(99,51){\tiny 4}
  \end{picture}
\end{center}
Suppose $b_4$ does not have $c_3$ as a favorite. Then we can
assume that $b_4$'s favorite is $c_4$ (remember that we know it is
not $c_1$). But then $a_1b_1c_1,a_2b_2c_2,a_3b_3c_3,a_4b_4c_4$ is
a stable matching since all $a_i$ and $b_i$ are 1-content, except
$a_2$ and $b_2$ which are already together. Thus we can assume
that $b_4$'s favorite is $c_3$.
\begin{center}
  \begin{picture}(130,80)(-4,-4)
    \sixteendots
    \put(0,0){\line(1,0){40}}
    \put(0,20){\line(1,0){80}}
    \put(0,40){\line(2,1){40}}
    \put(0,60){\line(1,0){120}}
    \put(80,60){\line(2,-1){40}}
    \put(40,0){\line(2,1){40}}
    \put(99,61){\tiny 3}
    \put(99,51){\tiny 4}
  \end{picture}
\end{center}
Again, form the matching
$a_1b_1c_1,a_2b_2c_2,a_3b_3c_3,a_4b_4c_4$. What are the possible
blocking triples? We observe that $a_2$ is the only $a_i$ which is
not 1-content, so a blocking triple must contain $a_2$.
Since $b_1$ and $b_3$ are 1-content and $a_2$ already has got $b_2$,
it follows that $b_4$ belongs to the blocking triple.
We also see that $c_1$ cannot belong to the blocking triple,
since $c_1$ already has $a_2$ which it prefers to $a_1$.
Thus, the only possible blocking triple
is $a_2b_4c_3$. In that case, $a_2$ must prefer $b_4$ to
$b_2$. In the same manner (using the matching
$a_1b_1c_1,a_2b_2c_2,a_3b_3c_4,a_4b_4c_3$ instead) we deduce that
$a_2$ prefers $b_3$ to $b_2$. This means that $a_2$'s second-best
choice is either $b_3$ or $b_4$. For symmetry reasons we can
assume it is $b_3$.
\begin{center}
  \begin{picture}(130,80)(-4,-4)
    \sixteendots
    \put(0,0){\line(1,0){40}}
    \put(0,20){\line(1,0){80}}
    \put(0,40){\line(2,1){40}}
    \put(0,60){\line(1,0){120}}
    \put(80,60){\line(2,-1){40}}
    \put(40,0){\line(2,1){40}}
    \dashline[-7]{3}(0,40)(40,20)
    \put(99,61){\tiny 3}
    \put(99,51){\tiny 4}
  \end{picture}
\end{center}
By the 113-condition we know that $c_3$'s favorite cannot be
$a_2$, $a_3$ or $a_4$, so it must be $a_1$.
\begin{center}
  \begin{picture}(130,80)(-4,-4)
    \sixteendots
    \put(0,0){\line(1,0){40}}
    \put(0,20){\line(1,0){80}}
    \put(0,40){\line(2,1){40}}
    \put(0,60){\line(1,0){120}}
    \put(80,60){\line(2,-1){40}}
    \put(40,0){\line(2,1){40}}
    \dashline[-7]{3}(0,40)(40,20)
    \put(80,20){\line(1,1){40}}
    \put(99,61){\tiny 3}
    \put(99,51){\tiny 4}
  \end{picture}
\end{center}
Using the 113-condition again, we see that $a_1$'s second-best
choice cannot be $b_3$ or $b_4$, so it must be $b_2$.
\begin{center}
  \begin{picture}(130,80)(-4,-4)
    \sixteendots
    \put(0,0){\line(1,0){40}}
    \put(0,20){\line(1,0){80}}
    \put(0,40){\line(2,1){40}}
    \put(0,60){\line(1,0){120}}
    \put(80,60){\line(2,-1){40}}
    \put(40,0){\line(2,1){40}}
    \dashline[-7]{3}(0,40)(40,20)
    \dashline[-7]{3}(0,60)(40,40)
    \put(99,61){\tiny 3}
    \put(99,51){\tiny 4}
  \end{picture}
\end{center}
Now, the matching $a_1b_2c_2,a_2b_1c_1,a_3b_3c_3,a_4b_4c_4$
is stable: The only $a_i$ who is
not 1-content is $a_1$. She wants to switch only to $b_1$ who is
1-content.

\section{Acknowledgment}\noindent
Research partially supported by the European Commission's IHRP
Programme, grant HPRN-CT-2001-00272, ``Algebraic Combinatorics in
Europe.''

\end{document}